\documentclass[11pt]{article}
\usepackage{mathrsfs}
\usepackage[centertags]{amsmath}
\usepackage{amsfonts}
\usepackage{amssymb}
\usepackage{amsthm}
\usepackage{cases}
\usepackage{indentfirst}
\usepackage{epsfig}
\usepackage{graphicx}
\usepackage{color}
\usepackage{array}
\allowdisplaybreaks[4]
\usepackage{setspace}
\date{}

\definecolor{c20}{rgb}{0.,0.7,0.}
\definecolor{c30}{rgb}{0.,0.,1.}
\definecolor{c40}{rgb}{1,0.1,0.7}
\definecolor{c50}{rgb}{1,0,0}
\definecolor{c60}{rgb}{0.00,0.00,1.00}

\def\lx#1{\textcolor{c30}{#1}}
\def\lx#1{#1}

\newtheorem{theorem}{Theorem}[section]

\newtheorem{lemma}{Lemma}[section]

\newtheorem{remark}{Remark}[section]

\numberwithin{equation}{section}

\def\P{\operatorname*{\mathbb{P}}}

\allowdisplaybreaks[4]

\begin{document}
\title{\textbf{ Second-order expansions for maxima of dynamic bivariate normal copulas}}
\author{$^{a}$Rui Wang  \qquad$^{b}$Xin Liao\thanks{Corresponding author. Email address: liaoxin2010@163.com} \qquad$^{a}$Zuoxiang Peng   \\
{\small $^{a}$School of Mathematics and Statistics, Southwest University, Chongqing, 400715, China}\\
 {\small $^{b}$Business School, University of Shanghai for Science and Technology, Shanghai, 200093, China}}

\maketitle
\begin{quote}
{\bf Abstract.}~~In this paper, we establish the second-order distributional expansions of normalized maxima of
$n$ independent observations, where the $i$th observation follows from a normal copula with its correlation coefficient being a monotone continuous function. These expansions can be used to deduce the convergence rates of
distributions of normalized maxima to their limits.

\end{quote}

\section{Introduction}
\label{sec1}

Let $\{(X_i,Y_i),1\leq i\leq n,n\geq 1\}$ denote independent and identically distributed bivariate random vectors with distribution function
$F(x,y)$ and continuous marginal distributions $F_1$ and $F_2$. The copula of $F$ is given by
$F(F_{1}^{-}(x),F_{2}^{-}(y))$, where $F_{i}^{-}$ denotes the inverse function of $F_{i}$, $i=1,2$.
We say that the copula of $F$ is a normal copula $C(x,y;\rho)$, if the density of $C(x,y;\rho)$ is given by
\begin{equation}\label{eq1.1}
c(x,y;\rho)=\frac{1}{\sqrt{1-\rho^2}} \exp\left( \frac{2\rho\Phi^{-}(x)\Phi^{-}(y)
-\rho^2(\Phi^{-}(x))^2-\rho^2(\Phi^{-}(y))^2}{2(1-\rho^2)} \right),
\end{equation}
where $\rho \in (-1,1)$ and $\Phi(x)$ is the standard normal distribution function.

Due to its easy to simulation and some attractive properties, the normal copula
has received many applications. Taylor et al. (2015) proposed causal quantities to evaluate surrogacy based on normal copula; Naldi and D'Acquisto (2008) considered the economic consequences of failures as a figure of merit of reliable communications networks by using normal copula, a few mentioned here.
But the biggest weakness of normal copula is its tail asymptotic independence,
see Sibuya (1960) and Embrechts et al. (2002). The tail asymptotic independence of normal copula may deduce
the under-estimation of extreme probabilities in risk management.
To overcome the drawback, Frick and Reiss (2013) showed that
\begin{eqnarray}\label{eq1.4}
&&\lim_{n\to \infty}\P\left(n(\max_{1\leq i \leq n}F_{1}(X_i)-1)\leq x,n(\max_{1\leq i \leq n}F_{2}(Y_i)-1)\leq y\right)\nonumber \\
&=& \exp\left(  \Phi\left( \sqrt{\lambda}+\frac{\log \frac{x}{y}}{2\sqrt{\lambda}} \right)x
+ \Phi\left( \sqrt{\lambda}+\frac{\log \frac{y}{x}}{2\sqrt{\lambda}} \right)y \right)
\end{eqnarray}
for $x<0$ and $y<0$, if the correlation coefficient
$\rho=\rho_n$ satisfies the following so-called H\"{u}sler-Reiss condition
\begin{eqnarray}\label{eq1.3}
(1-\rho_n)\log n \to \lambda \in [0,\infty]
\end{eqnarray}
as $n\to \infty$, see H\"{u}sler and Reiss (1989). Note that \eqref{eq1.4} is a copula version of the limit in H\"{u}sler and Reiss (1989) for the normalized
maxima of bivariate normal triangular arrays with correlation coefficients satisfying \eqref{eq1.3}.
Further, Liao et al. (2016) extended the work of Frick and Reiss (2013) by assuming
\begin{eqnarray}\label{eq1.2}
  \rho_{ni}=1-\frac{m\left( i/n \right)}{\log n}
\end{eqnarray}
for some nonnegative function $m(x)$, which allows $\rho_{ni}$ depending on both $i$ and $n$. Under the condition
\eqref{eq1.2}, Liao et al. (2016) proved that
\begin{eqnarray}\label{eq1.5}
\lim_{n\to \infty} \P\left(n\left(\max_{1\leq i \leq n}F_1(X_{i})-1\right)\leq x,n\left(\max_{1\leq i \leq n}F_2(Y_{i})-1\right)\leq y\right)
=G(x,y),
\end{eqnarray}
with
\begin{eqnarray*}
G(x,y)=\exp\left( x\int_{0}^{1} \Phi\left( \sqrt{m(s)} +\frac{\log \frac{x}{y}}{2\sqrt{m(s)}} \right)ds
+ y\int_{0}^{1} \Phi\left( \sqrt{m(s)} +\frac{\log \frac{y}{x}}{2\sqrt{m(s)}} \right)ds  \right)
\end{eqnarray*}
if $m(s)$ defined on $[0,1]$ is continuous and positive, and $G(x,y)=\exp(x+y)$ as $\lim_{n\to \infty}\min_{1\leq i \leq n}m\left( i/n \right)=\infty$, and $G(x,y)=\exp(\min(x,y))$ if $\lim_{n\to \infty}\max_{1\leq i \leq n} m\left( i/n \right)= 0$. The above normal copulas with
$\rho_{ni}$ depending on both $i$ and $n$ or only sample size $n$, can be called dynamic copulas.
Recently, various dynamic copulas are receiving more attention in modeling financial time series; see, e.g., Salvatierra and Patton (2015),
Wang et al. (2015), Chu (2015).

In this paper, we are interested in the second-order distributional expansions of normalized maxima
$(n\left(\max_{1\leq i \leq n}F_1(X_{i})-1\right)\leq x,n\left(\max_{1\leq i \leq n}F_2(Y_{i})-1\right)\leq y)$.
For independent bivariate normal triangular arrays satisfying \eqref{eq1.3}, Hashorva et al. (2016) imposed second-order H\"{u}sler-Reiss
condition, and derived the second-order distributional expansions of maxima. Under the second-order
H\"{u}sler-Reiss condition, Liao and Peng (2014, 2015) obtained the uniform convergence rates of maxima,
and the second-order expansions of joint distributions of maxima and minima of bivariate normal triangular arrays.
For the independent and non-identically distributed bivariate normal triangular arrays satisfying \eqref{eq1.2},
the second-order distributional expansions of maxima are given by Liao and Peng (2016), and the second-order
expansions of joint distributions of maxima and minima are derived by Lu and Peng (2017).
To the best of our knowledge, there are no studies on the second-order expansions of distributions of
$(n(\max_{1\leq i \leq n}F_1(X_i)-1)\leq x,n(\max_{1\leq i \leq n}F_2(Y_i)-1)\leq y)$. The aim of this
paper is to fill this gap.

The rest of this paper is organized as follows. In Section \ref{sec2}, we establish the second-order
distributional expansions of $(n(\max_{1\leq i \leq n}F_1(X_i)-1)\leq x,n(\max_{1\leq i \leq n}F_2(Y_i)-1)\leq y)$
by considering three cases: m(s) is continuous positive function on [0,1], $\lim_{n\to \infty}\min_{1\leq i \leq n}m\left( i/n \right)=\infty$, and $\lim_{n\to \infty}\max_{1\leq i \leq n} m\left( i/n \right)= 0$. Some examples are also given in Section \ref{sec2}.  All proofs are deferred in Section \ref{sec3}.

\section{Main Results}
\label{sec2}

In this section, the second-order expansions of
\[G_{n}(x,y)=\P\Big(n(\max_{1\leq i \leq n}F_{1}(X_i)-1)\leq x,
n(\max_{1\leq i \leq n}F_{2}(Y_i)-1)\leq y \Big)\]
are established with $m(s)$ given by \eqref{eq1.2} satisfying some regular conditions. Note that Liao et al. (2016) derived the convergence of $G_{n}(x,y)$ for three different cases.
In order to derive the second-order expansions of $G_{n}(x,y)$, additional conditions are needed in each case.
The following Theorem is about the second-order distributional expansion of $G_{n}(x,y)$ as $m(s)$ defined on $[0,1]$ is monotone, continuous and positive.

\begin{theorem}
\label{the2.1}
Assume that \eqref{eq1.2} holds with $m(s)$ being monotone and continuous on $[0,1]$. For any $x<0$ and $y<0$, we have
\begin{equation}
\label{eq2.1}
\lim_{n\rightarrow\infty} \frac{\log n}{\log \log n}\Big[G_{n}(x,y)-G(x,y)\Big]
=\frac{1}{2\sqrt{2\pi}} G(x,y) \int_0^1 \sqrt{m(s)} \exp\left(-\frac{m(s) \lx{-} \log(xy)+\frac{(\log({x}/{y}))^2}{4m(s)}}{2}\right) ds.
\end{equation}

\end{theorem}

\noindent{\bf Example 2.1.}~
Assume that \eqref{eq1.2} holds with $m(x)=x+1$, $0\leq x\leq 1$. It follows from Theorem \ref{the2.1}
that
\begin{eqnarray*}
G_{n}(x,y)=G(x,y)+\frac{\log \log n}{\log n}\lx{\frac{(-x)G(x,y)}{\sqrt{2\pi}}}\int_{1}^{\sqrt{2}}s^2
\exp\left( -\frac{\left( s+\frac{\log ({x}/{y})}{2s} \right)^2}{2} \right)ds (1+o(1))
\end{eqnarray*}
for large $n$.

For the case of $\lim_{n\to \infty} \min_{1\leq i \leq n}m(i/n)=\infty$, with addition condition Theorem \ref{the2.2} shows the second-order expansion of $G_{n}(x,y)$ as follows.

\begin{theorem}
\label{the2.2}
Assume that \eqref{eq1.2} holds with $\lim_{n\to \infty} \min_{1\leq i \leq n}m(i/n)=\infty$
and $\lim_{n\rightarrow\infty}\frac{\max_{1\leq i \leq n} m(i/n)}{\log\log n}=0$, then
for any $x<0$ and $y<0$ we have
\begin{equation}
\label{eq2.2}
\lim_{n\to\infty}\left( \frac{1}{n} \sum_{i=1}^n \frac{ \exp\left( -{m(i/n)}/{2} \right)}{\sqrt{m(i/n)}} \right)^{-1}\Big[G_{n}(x,y)-e^{x+y}\Big]
=
\sqrt{\frac{2}{\pi}} (xy)^{\frac{1}{2}} e^{x+y}.
\end{equation}
\end{theorem}

\noindent{\bf Example 2.2.}~
Let
\begin{equation*}
m(i/n)=\left\{\begin{array}{lll}
&4\log \log \log \frac{n}{i},& i \in [1,n^{\frac{1}{2}}],\\
&2\left( \log \log \log n \bigvee \log \log \log \frac{n}{i} \right),& \text{otherwise}.
\end{array}\right.\end{equation*}
One can check that $m(i/n)$ given above satisfies the conditions of Theorem \ref{the2.2}, so by Theorem \ref{the2.2}, the second-order expansion of
$G_{n}(x,y)$ is given by
\begin{equation*}
G_{n}(x,y)=e^{x+y}+\frac{\sqrt{\frac{2}{\pi}}(xy)^{\frac{1}{2}}e^{x+y}}{(\log \log n)\sqrt{2\log \log \log n}}(1+o(1))
\end{equation*}
for large $n$.

For the last case of $\lim_{n\to \infty}\max_{1 \leq i \leq n} m(i/n)=0$, with additional condition we have the following second-order expansion of $G_{n}(x,y)$.

\begin{theorem}
\label{the2.3}
Assume that \eqref{eq1.2} holds with $\lim_{n\to \infty} (\log \log n )\min_{1\leq i \leq n}m(i/n)=\infty$ and
\\$\lim_{n\to \infty}\max_{1 \leq i \leq n} m(i/n)=0$. Then for $x<0, y<0$,
\begin{itemize}
\item[(i)]~ if $x\neq y$, we have
\begin{equation}
\label{eq2.3}
\lim_{n\to\infty}\left( \frac{1}{n} \sum_{i=1}^{n} \left(m(i/n)\right)^{\frac{3}{2}}  \exp\left( -\frac{\left( \log \frac{\min(x,y)}{\max(x,y)}\right)^2}{8m(i/n)} \right)  \right)^{-1}\Big[G_{n}(x,y)-e^{\min(x,y)}\Big]=
-\frac{8 (xy)^{\frac{1}{2}}e^{\min(x,y)}}{\sqrt{2\pi} \log \frac{\min(x,y)}{\max(x,y)}}.
\end{equation}

\item[(ii)]~if $x=y$, we have
\begin{equation}
\lim_{n\to \infty}\left( \frac{1}{n}\sum_{i=1}^{n}\sqrt{m(i/n)} \right)^{-1}
\left[ G_n(x,y)-e^{x} \right]=-\frac{2x}{\pi}e^{x}.
\end{equation}
\end{itemize}
\end{theorem}

\noindent{\bf Example 2.3.}~
One can check that
\begin{equation*}
m(i/n)=\left\{\begin{array}{lll}
&\frac{(\log \log i)^{\frac{1}{2}}}{\log \log \log n},& i \in [1,\log n],\\
&\frac{1}{\log \log \log n}\bigwedge\frac{(\log \log i)^{\frac{1}{2}}}{\log \log \log n} ,& \text{otherwise}
\end{array}\right.\end{equation*}
satisfies the conditions of Theorem \ref{the2.3}, and the second-order expansion of $G_n(x,y)$
is given by
\begin{equation*}
G_{n}(x,y)=\left\{\begin{array}{lll}
&e^{\min(x,y)}-\frac{1}{(\log \log \log n)^{\frac{3}{2}}(\log \log n)^{\frac{1}{8}\left(\log \frac{\min(x,y)}{\max(x,y)}\right)^2}}\cdot
\frac{8(xy)^{\frac{1}{2}}e^{\min(x,y)}}{\sqrt{2\pi}\log \frac{\min(x,y)}{\max(x,y)}}(1+o(1)),& x\neq y,\\
&e^{x}-\frac{2xe^{x}}{\pi (\log \log \log n)^{\frac{1}{2}}}(1+o(1)) ,& x=y
\end{array}\right.\end{equation*}
for large $n$.

\begin{remark}
For different cases, Theorems \ref{the2.1}-\ref{the2.3} show that the convergence rates of $G_n(x,y)$
to $G(x,y)$ are given as follows:
\begin{itemize}
\item[(i)]~if $m(s)$ is monotone and continuous, Theorem \ref{the2.1} shows that
the convergence rate is proportional to
$\frac{\log \log n}{\log n}$.

\item[(ii)]~if $m(s)$ satisfies $\lim_{n\to \infty} \min_{1\leq i \leq n}m(i/n)=\infty$
and $\lim_{n\rightarrow\infty}\frac{\max_{1\leq i \leq n} m(i/n)}{\log\log n}=0$, Theorem \ref{the2.2} shows that
the convergence rate is the same order of
$\frac{1}{n}\sum_{i=1}^{n}\frac{\exp(-m(i/n)/2)}{\sqrt{m(i/n)}} $.

\item[(iii)]~if $m(s)$ satisfies $\lim_{n\to \infty}\max_{1 \leq i \leq n} m(i/n)=0$ and $\lim_{n\to \infty} (\log \log n )\min_{1\leq i \leq n}m(i/n)=\infty$, Theorem \ref{the2.3} shows that
the convergence rate of $G_n(x,y)$ to its limit $G(x,y)$ is the same order of $\frac{1}{n}\sum_{i=1}^{n} (m(i/n))^{\frac{3}{2}}\exp\left( -\frac{\log \frac{\min(x,y)}{\max(x,y)}}{8m(i/n)} \right)$ for $x\neq y$, and the same order
of $\frac{1}{n}\sum_{i=1}^{n}\sqrt{m(i/n)}$ for $x=y$.
\end{itemize}
\end{remark}

\section{Proofs}
\label{sec3}
\noindent

The aim of this section is to prove our main results. In order to prove Theorem \ref{the2.1},
we need the following key lemma, which shows the convergence rate of
$\frac{1}{n}\sum_{i=1}^n\int_{y}^{0}\Phi\left( \sqrt{m(i/n)}+\frac{\log (t/x)}{2\sqrt{m(i/n)}} \right)dt$.

\begin{lemma}
\label{lem3.1}
Under the conditions of Theorem \ref{the2.1}, for $x<0$ and $y<0$ we have
\begin{eqnarray}
\label{eq3.1}
\int_0^1 \int_y^0
\Phi\left(\sqrt{m(s)}+\frac{\log(t/x)}{2\sqrt{m(s)}}\right)dtds
-\frac{1}{n}\sum_{i=1}^{n}\int_y^0\Phi\left(\sqrt{m(i/n)}+\frac{\log(t/x)}{2\sqrt{m(i/n)}}\right)dt
=O\left(\frac{1}{n}\right).
\end{eqnarray}
\end{lemma}

\noindent
{\bf Proof of Lemma \ref{lem3.1}.}~~
Without loss of generality, assume that $m(s)$ is increasing.

For $x\leq y$, noting that $\int_y^{0}\Phi\Big(\sqrt{m(s)}+\frac{\log(t/x)}{2\sqrt{m(s)}}\Big)dt$ is increasing about $s$, we have
\begin{eqnarray}
\label{eq3.2}
\int_0^1\int_y^0\Phi\left(\sqrt{m(s)}+\frac{\log(t/x)}{2\sqrt{m(s)}}\right)dtds
&=&\sum_{i=1}^n \int_{\frac{i-1}{n}}^{\frac{i}{n}}
\int_y^0\Phi\left(\sqrt{m(s)}+\frac{\log(t/x)}{2\sqrt{m(s)}}\right)dtds\nonumber
\\
&<&\sum_{i=1}^n \int_{\frac{i-1}{n}}^{\frac{i}{n}}
\int_y^0\Phi\left(\sqrt{m(i/n)}+\frac{\log(t/x)}{2\sqrt{m(i/n)}}\right)dtds\nonumber
\\
&=&\frac{1}{n}\sum_{i=1}^n\int_y^0\Phi\left(\sqrt{m(i/n)}+\frac{\log(t/x)}{2\sqrt{m(i/n)}}\right)dt
\end{eqnarray}
and
\begin{align}
\label{eq3.3}
\int_0^1\int_y^0\Phi\left(\sqrt{m(s)}+\frac{\log(t/x)}{2\sqrt{m(s)}}\right)dtds
&=\sum_{i=0}^{n-1} \int_{\frac{i}{n}}^{\frac{i+1}{n}}
\int_y^0\Phi\left(\sqrt{m(s)}+\frac{\log(t/x)}{2\sqrt{m(s)}}\right)dtds\nonumber
\\
&>\sum_{i=0}^{n-1} \int_{\frac{i}{n}}^{\frac{i+1}{n}}
\int_y^0\Phi\left(\sqrt{m(i/n)}+\frac{\log(t/x)}{2\sqrt{m(i/n)}}\right)dtds\nonumber
\\
&=\frac{1}{n}\sum_{i=1}^{n}\int_y^0\Phi\left(\sqrt{m(i/n)}+\frac{\log(t/x)}{2\sqrt{m(i/n)}}\right)dt
+O\left(\frac{1}{n}\right),
\end{align}
so, \eqref{eq3.2} and \eqref{eq3.3} implies that \eqref{eq3.1} holds for $x\leq y$.

For $x>y$, let
\begin{eqnarray*}
\int_0^1 \int_y^0
\Phi\left(\sqrt{m(s)}+\frac{\log(t/x)}{2\sqrt{m(s)}}\right)dtds
-\frac{1}{n} \sum_{i=1}^{n} \int_y^0
\Phi\left(\sqrt{m(i/n)}+\frac{\log(t/x)}{2\sqrt{m(i/n)}}\right)dt
=
A_1(n)+A_2(n)
\end{eqnarray*}
with
\begin{eqnarray*}
A_1(n)=\int_0^1 \int_x^0
\Phi\left(\sqrt{m(s)}+\frac{\log(t/x)}{2\sqrt{m(s)}}\right)dtds
-\frac{1}{n} \sum_{i=1}^{n} \int_x^0
\Phi\left(\sqrt{m(i/n)}+\frac{\log(t/x)}{2\sqrt{m(i/n)}}\right)dt
\end{eqnarray*}
and
\begin{eqnarray}\label{eq3.4}
A_2(n)=\int_0^1 \int_y^x
\Phi\left(\sqrt{m(s)}+\frac{\log(t/x)}{2\sqrt{m(s)}}\right)dtds
-\frac{1}{n} \sum_{i=1}^{n} \int_y^x
\Phi\left(\sqrt{m(i/n)}+\frac{\log(t/x)}{2\sqrt{m(i/n)}}\right)dt.
\end{eqnarray}
By arguments similar to \eqref{eq3.2} and \eqref{eq3.3}, we can get
$
A_1(n)=O({1}/{n}).
$
The rest is to show that
$
A_2(n)=O({1}/{n}).
$

First note that the function $f(z)=\sqrt{z}+\frac{\log(y/x)}{2\sqrt{z}}$  is increasing  for $z>\frac{\log(y/x)}{2}$ and decreasing as $z<\frac{\log(y/x)}{2}$. So, we need to deal with \eqref{eq3.4}
through the following three cases:
\[
\mbox{(i).}\quad y\leq xe^{2m(1)};\qquad
\mbox{(ii).}\quad xe^{2m(1)}<y<xe^{2m(0)};\qquad
\mbox{(iii).}\quad xe^{2m(0)}\leq y\leq x.\]
Arguments similar to that of \eqref{eq3.2} and \eqref{eq3.3}, we can show that $A_{2}(n)=O(1/n)$ for case (iii). Details are omitted here. So, there are only cases (i) and (ii) left as we estimate the bound of \eqref{eq3.4}.

For case (i), note that $y\leq xe^{m(s)}\leq x<0$ for $s\in [0,1]$ since $m(s)$ is increasing and
$\Phi\Big(\sqrt{m(s)}+\frac{\log(t/x)}{2\sqrt{m(s)}}\Big)$ is decreasing respect to $s$ for $t\in [y,xe^{2m(s)}]$. Hence,
\begin{align}
\label{eq3.5}
\int_0^1\int_y^{xe^{2m(s)}}\Phi\left(\sqrt{m(s)}+\frac{\log(t/x)}{2\sqrt{m(s)}}\right)dtds
&=\sum_{i=1}^n \int_{\frac{i-1}{n}}^{\frac{i}{n}}
\int_y^{xe^{2m(s)}}\Phi\left(\sqrt{m(s)}+\frac{\log(t/x)}{2\sqrt{m(s)}}\right)dtds\nonumber
\\
&>\sum_{i=1}^n \int_{\frac{i-1}{n}}^{\frac{i}{n}}
\int_y^{xe^{2m(i/n)}}\Phi\left(\sqrt{m(i/n)}+\frac{\log(t/x)}{2\sqrt{m(i/n)}}\right)dtds\nonumber
\\
&=\frac{1}{n}\sum_{i=1}^n
\int_y^{xe^{2m(i/n)}}\Phi\left(\sqrt{m(i/n)}+\frac{\log(t/x)}{2\sqrt{m(i/n)}}\right)dt
\end{align}
and
\begin{align}
\label{eq3.6}
\int_0^1\int_y^{xe^{2m(s)}}\Phi\left(\sqrt{m(s)}+\frac{\log(t/x)}{2\sqrt{m(s)}}\right)dtds
&=\sum_{i=0}^{n-1} \int_{\frac{i}{n}}^{\frac{i+1}{n}}
\int_y^{xe^{2m(s)}}\Phi\left(\sqrt{m(s)}+\frac{\log(t/x)}{2\sqrt{m(s)}}\right)dtds
\nonumber
\\
&<\sum_{i=0}^{n-1} \int_{\frac{i}{n}}^{\frac{i+1}{n}}
\int_y^{xe^{2m(s)}}\Phi\left(\sqrt{m(i/n)}+\frac{\log(t/x)}{2\sqrt{m(i/n)}}\right)dtds\nonumber
\\
&<\sum_{i=0}^{n-1} \int_{\frac{i}{n}}^{\frac{i+1}{n}}
\int_y^{xe^{2m(i/n)}}\Phi\left(\sqrt{m(i/n)}+\frac{\log(t/x)}{2\sqrt{m(i/n)}}\right)dtds
\nonumber
\\
&=\frac{1}{n}\sum_{i=1}^{n}\int_y^{xe^{2m(i/n)}}\Phi\left(\sqrt{m(i/n)}+\frac{\log(t/x)}{2\sqrt{m(i/n)}}\right)dt
+O\left(\frac{1}{n}\right).
\end{align}
Combining \eqref{eq3.5} with  \eqref{eq3.6}, we have
\begin{equation}
\label{eq3.7}
\int_0^1\int_y^{xe^{2m(s)}}\Phi\left(\sqrt{m(s)}+\frac{\log(t/x)}{2\sqrt{m(s)}}\right)dtds
-
\frac{1}{n}\sum_{i=1}^n
\int_y^{xe^{2m(i/n)}}\Phi\left(\sqrt{m(i/n)}+\frac{\log(t/x)}{2\sqrt{m(i/n)}}\right)dt
=
O\left(\frac{1}{n}\right).
\end{equation}

Similarly, noting that $\Phi\left(\sqrt{m(s)}+\frac{\log(t/x)}{2\sqrt{m(s)}}\right)$ is increasing respect to $s$ for $t\in (xe^{2m(s)},x]$, we have
\begin{eqnarray*}
\int_0^1\int_{xe^{2m(s)}}^x\Phi\left(\sqrt{m(s)}+\frac{\log(t/x)}{2\sqrt{m(s)}}\right)dtds
<
\frac{1}{n}\sum_{i=1}^n
\int_{xe^{2m(i/n)}}^x\Phi\left(\sqrt{m(i/n)}+\frac{\log(t/x)}{2\sqrt{m(i/n)}}\right)dt
\end{eqnarray*}
and
\begin{equation*}
\int_0^1\int_{xe^{2m(s)}}^x\Phi\left(\sqrt{m(s)}+\frac{\log(t/x)}{2\sqrt{m(s)}}\right)dtds
>
\frac{1}{n}\sum_{i=1}^n
\int_{xe^{2m(i/n)}}^x\Phi\left(\sqrt{m(i/n)}+\frac{\log(t/x)}{2\sqrt{m(i/n)}}\right)dt
+O\left(\frac{1}{n}\right),
\end{equation*}
implying
\begin{equation}
\label{eq3.8}
\int_0^1\int_{xe^{2m(s)}}^x\Phi\left(\sqrt{m(s)}+\frac{\log(t/x)}{2\sqrt{m(s)}}\right)dtds
-
\frac{1}{n}\sum_{i=1}^n
\int_{xe^{2m(i/n)}}^x\Phi\left(\sqrt{m(i/n)}+\frac{\log(t/x)}{2\sqrt{m(i/n)}}\right)dt
=
O\left(\frac{1}{n}\right).
\end{equation}
It follows from \eqref{eq3.7} and \eqref{eq3.8} that \eqref{eq3.4} holds for case (i).

For case (ii), i.e. $xe^{2m(1)}<y<xe^{2m(0)}$, there exists $s_0\in (0,1)$ such
that $y=xe^{2m(s_0)}$ since $m(s)$ is increasing  and continuous. We split the following integral into two parts:
\begin{eqnarray}
\label{eq3.9}
& & \int_0^1 \int_y^{xe^{2m(s)}} \Phi\left( \sqrt{m(s)} +\frac{\log(t/x)}{2\sqrt{m(s)}}  \right) dt ds \nonumber \\
&=& \int_0^{s_0} \int_y^{xe^{2m(s)}} \Phi\left( \sqrt{m(s)} +\frac{\log(t/x)}{2\sqrt{m(s)}}  \right) dt ds
+\int_{s_0}^1 \int_y^{xe^{2m(s)}} \Phi\left( \sqrt{m(s)} +\frac{\log(t/x)}{2\sqrt{m(s)}}  \right) dt ds. \nonumber \\
\end{eqnarray}

By arguments similar to \eqref{eq3.5}-\eqref{eq3.8}, we have
\begin{eqnarray*}
\int_0^{s_0} \int_y^{xe^{2m(s)}} \Phi\left( \sqrt{m(s)} +\frac{\log(t/x)}{2\sqrt{m(s)}}  \right) dt ds
= \frac{1}{n}\sum_{i=1}^{[ns_0]}\int_y^{xe^{2m(i/n)}} \Phi\left( \sqrt{m(i/n)}
+ \frac{\log(t/x)}{2\sqrt{m(i/n)}}  \right) +O\left(\frac{1}{n}\right)
\end{eqnarray*}
and
\begin{eqnarray*}
\int_{s_0}^{1} \int_y^{xe^{2m(s)}} \Phi\left( \sqrt{m(s)} +\frac{\log(t/x)}{2\sqrt{m(s)}}  \right) dt ds
= \frac{1}{n}\sum_{i=[ns_0]+1}^{n}\int_y^{xe^{2m(i/n)}} \Phi\left( \sqrt{m(i/n)}
+ \frac{\log(t/x)}{2\sqrt{m(i/n)}}  \right) +O\left(\frac{1}{n}\right),
\end{eqnarray*}
implying that
\begin{equation}
\label{eq3.10}
\int_0^1 \int_y^{xe^{2m(s)}} \Phi\left( \sqrt{m(s)} +\frac{\log(t/x)}{2\sqrt{m(s)}}  \right) dt ds
= \frac{1}{n}\sum_{i=1}^{n}\int_y^{xe^{2m(i/n)}} \Phi\left( \sqrt{m(i/n)}
+ \frac{\log(t/x)}{2\sqrt{m(i/n)}}  \right) +O\left(\frac{1}{n}\right).
\end{equation}
It follows from \eqref{eq3.8} and \eqref{eq3.10} that \eqref{eq3.4} holds for case (ii).

The proof is complete.
\qed

In order to prove Theorems \ref{the2.1}-\ref{the2.3}, we first give the following definitions:
\begin{eqnarray}
\label{eq3.11}
I_k(x,y;m(i/n))=\int_y^{-\frac{1}{\log n}}\left(-\log(-t)\right)^k
\phi\left(\sqrt{m(i/n)}+\frac{\log(t/x)}{2\sqrt{m(i/n)}}\right)dt,
k=0,1,2,3,
\end{eqnarray}
where $\phi(x)$ is the standard normal density.
One can check that
\begin{equation}
\label{eq3.12}
I_0(x,y;m(i/n))=-2x\sqrt{m(i/n)}\left[\Phi\left(\sqrt{m(i/n)}+\frac{\log(-x\log{n})}{2\sqrt{m(i/n)}}\right)
-\Phi\left(\sqrt{m(i/n)}+\frac{\log(x/y)}{2\sqrt{m(i/n)}}\right)\right],
\end{equation}
\begin{eqnarray}
\label{eq3.13}
I_1(x,y;m(i/n))&=&4xm(i/n)\left[\phi\left(\sqrt{m(i/n)}+\frac{\log(-x\log{n})}{2\sqrt{m(i/n)}}\right)
-\phi\left(\sqrt{m(i/n)}+\frac{\log(x/y)}{2\sqrt{m(i/n)}}\right)\right]\nonumber\\
&&+2x\sqrt{m(i/n)}\Big(\log(-x)+2m(i/n)\Big)
\left[\Phi\left(\sqrt{m(i/n)}+\frac{\log(-x\log{n})}{2\sqrt{m(i/n)}}\right) \right. \nonumber\\
&&
\left.
-\Phi\left(\sqrt{m(i/n)}+\frac{\log(x/y)}{2\sqrt{m(i/n)}}\right)\right],
\end{eqnarray}
\begin{eqnarray}
\label{eq3.14}
I_2(x,y;m(i/n))&=&-2x\sqrt{m(i/n)}\left(4(m(i/n))^2+4(m(i/n))+4(m(i/n))\log(-x)+(\log(-x))^2\right)\nonumber\\
&&
\times\left[\Phi\left(\sqrt{m(i/n)}+\frac{\log(-x\log{n})}{2\sqrt{m(i/n)}}\right)
-\Phi\left(\sqrt{m(i/n)}+\frac{\log(x/y)}{2\sqrt{m(i/n)}}\right)\right]\nonumber\\
&&
-4x \left( 2(m(i/n))^2 +m(i/n)\log (-x/\log n) \right)
\phi\left( \sqrt{m(i/n)}+\frac{\log (-x\log n)}{2\sqrt{m(i/n)}} \right) \nonumber \\
&&
+4x \left( 2(m(i/n))^2 +m(i/n)\log (xy) \right)  \phi\left(\sqrt{m(i/n)}
+\frac{\log (x/y)}{2\sqrt{m(i/n)}} \right),
\end{eqnarray}
\begin{eqnarray}
\label{eq3.15}
I_3(x,y;m(i/n))&=& 2x\sqrt{m(i/n)}\Big( 8(m(i/n))^3 +(24+12\log (-x))(m(i/n))^2  \nonumber\\
& &
+ \left(  6(\log (-x))^2+12\log (-x) \right)m(i/n) +(\log (-x))^3 \Big) \nonumber\\
& &
\times \left[ \Phi\left( \sqrt{m(i/n)} +\frac{\log (-x\log n)}{2\sqrt{m(i/n)}} \right) -\Phi\left( \sqrt{m(i/n)} +\frac{\log (\frac{x}{y})}{2\sqrt{m(i/n)}}  \right) \right] \nonumber\\
& & +2x \left(  8(m(i/n))^3 +(m(i/n))^2(8\log (-x) -4\log \log n +16) \right.  \nonumber\\
&&\left.\ \ +2 m(i/n)((\log (-x))^2 -(\log \log n)\log (-x) +(\log \log n)^2 )\right)\nonumber\\
&& \times\phi\left( \sqrt{m(i/n)} +\frac{\log (-x\log n)}{2\sqrt{m(i/n)}} \right)
\nonumber\\
&&-2x \left(  8(m(i/n))^3 +(m(i/n))^2 (8\log (-x) +4\log(-y) +16) \right. \nonumber\\
&&\left.\ \ +2m(i/n)((\log (-x))^2 +(\log (-y))\log (-x) +(\log (-y))^2 )\right)\nonumber\\
&&\times \phi\left( \sqrt{m(i/n)} +\frac{\log (\frac{x}{y})}{2\sqrt{m(i/n)}} \right).
\end{eqnarray}
With Lemma 3.1, we can prove Theorem \ref{the2.1} as follows.

\noindent
{\bf Proof of Theorem \ref{the2.1}.}~~
By the Mill's ratio of normal distribution, for any fixed $x<0$ we have,
\begin{eqnarray}
\label{eq3.16}
\Phi^{-}(1+\frac{x}{n})&=&\sqrt{2\log{n}}
\left( 1-\frac{\log 4\pi +\log \log n}{4\log n} +\frac{\log 4\pi +\log \log n}{8(\log n)^2}
-\frac{(\log 4\pi +\log \log n)^2}{32(\log n)^2}  \right)\nonumber
\\
&&
-\frac{\log(-x)}{\sqrt{2\log{n}}}
\left( 1-\frac{1}{2\log n} +\frac{\log (-x)}{4\log n} +\frac{\log 4\pi +\log \log n}{4\log n} \right)
+o\Big({(\log{n})^{-\frac{3}{2}}}\Big).
\end{eqnarray}
Note that $o\Big({(\log{n})^{-\frac{3}{2}}}\Big)$ also holds uniformly for $x\in[y,-\frac{1}{\log n}]$ with fixed $y$.
It follows from \eqref{eq1.2}, \eqref{eq3.16} and the monotonicity and continuity of $m(s)$ that for large $n$ and fixed $x<0$ and $y<0$,
\begin{eqnarray}
\label{eq3.18}
&&
\frac{\Phi^{-}(1+\frac{x}{n})-\rho_{ni}\Phi^{-}(1+\frac{t}{n})}{\sqrt{1-\rho_{ni}^2}}\nonumber
\\
&=&\frac{\Phi^{-}(1+\frac{x}{n})-\Phi^{-}(1+\frac{t}{n})}{\sqrt{1-\rho_{ni}^2}}
+ \sqrt{\frac{1-\rho_{ni}}{1+\rho_{ni}}}\Phi^{-}(1+\frac{t}{n})\nonumber\\
&=&\sqrt{m(i/n)}+\frac{\log(t/x)}{2\sqrt{m(i/n)}}
- \frac{\log \log n}{4\log n}\left(  \sqrt{m(i/n)} - \frac{\log (t/x)}{2\sqrt{m(i/n)}}\right) \nonumber\\
& & +\frac{(\log (-t))^2 +\left( \log 4\pi -3m(i/n) -2 \right)\log (-t)}{8\sqrt{m(i/n)}\log n}
- \frac{(\log (-x))^2 +\left( \log 4\pi +m(i/n) -2 \right)\log (-x)}{8\sqrt{m(i/n)}\log n} \nonumber\\
& & +\frac{(m(i/n))^{\frac{3}{2}}-(\log 4\pi)\sqrt{m(i/n)}}{4\log n} +\frac{1+\log (-t)}{\log n} o(1)
\end{eqnarray}
holds uniformly for all $1\leq i \leq n$ and $t \in [y,-\frac{1}{\log n}]$.
Noting that
\begin{eqnarray}
\label{eq3.19}
&&
\int_{y}^{-\frac{1}{\log n}}
\phi\Big(\sqrt{m(i/n)}+\frac{\log(t/x)}{2\sqrt{m(i/n)}}\Big)
\Big(\sqrt{m(i/n)}-\frac{\log(t/x)}{2\sqrt{m(i/n)}}\Big)dt\nonumber
\\
&=&\int_{y}^{-\frac{1}{\log n}}
\frac{x}{t}\phi\Big(\sqrt{m(i/n)}-\frac{\log(t/x)}{2\sqrt{m(i/n)}}\Big)
\Big(\sqrt{m(i/n)}-\frac{\log(t/x)}{2\sqrt{m(i/n)}}\Big)dt\nonumber
\\
&=&2\sqrt{m(i/n)}x\int_{y}^{-\frac{1}{\log n}}
d\Big(\phi\Big(\sqrt{m(i/n)}-\frac{\log(t/x)}{2\sqrt{m(i/n)}}\Big)\Big)\nonumber
\\
&=&2\sqrt{m(i/n)}x
\Big[\phi\Big(\sqrt{m(i/n)}+\frac{\log(-x\log{n})}{2\sqrt{m(i/n)}}\Big)
-\phi\Big(\sqrt{m(i/n)}+\frac{\log(x/y)}{2\sqrt{m(i/n)}}\Big)\Big]\nonumber\\
&=& 2(-x)\sqrt{m(i/n)}\phi\left(  \sqrt{m(i/n)} +\frac{\log (x/y)}{2\sqrt{m(i/n)}} \right)
+ O\left({(\log n)^{-\frac{1}{2}}} \right) \nonumber \\
&=&
\sqrt{\frac{2}{\pi}}\sqrt{m(i/n)}\exp\left(-\frac{m(i/n)\lx{-} \log(xy)+\frac{(\log (x/y))^2}{4m(i/n)}}{2} \right)
+ O\left({(\log n)^{-\frac{1}{2}}} \right),
\end{eqnarray}
we have
\begin{eqnarray}
\label{eq3.20}
&&
-\frac{1}{n}\sum_{i=1}^{n}\int_{y}^{-\frac{1}{\log n}}
\phi\left(\sqrt{m(i/n)}+\frac{\log(t/x)}{2\sqrt{m(i/n)}}\right)
\left(\frac{\Phi^{-}(1+\frac{x}{n})-\rho_{ni}\Phi^{-}(1+\frac{t}{n})}{\sqrt{1-\rho_{ni}^2}}-\sqrt{m(i/n)}
-\frac{\log(t/x)}{2\sqrt{m(i/n)}}\right)dt\nonumber
\\
&=&\frac{\log\log{n}}{4\log{n}}
\frac{1}{n}\sum_{i=1}^{n}\int_{y}^{-\frac{1}{\log n}}
\phi\left(\sqrt{m(i/n)}+\frac{\log(t/x)}{2\sqrt{m(i/n)}}\right)
\left(\sqrt{m(i/n)}-\frac{\log(t/x)}{2\sqrt{m(i/n)}}
\right)dt\nonumber
\\
&&-\frac{1}{4n\log n}\sum_{i=1}^{n}\frac{1}{2\sqrt{m(i/n)}}I_2(x,y;m(i/n))
+\frac{1}{4n \log n} \sum_{i=1}^{n}\frac{\log 4\pi -3m(i/n)-2}{2\sqrt{m(i/n)}}I_1(x,y;m(i/n))\nonumber
\\
&&+\frac{1}{4n\log n}\nonumber\\
&&\quad \times\sum_{i=1}^{n} \frac{(\log (-x))^2+(\log 4\pi+m(i/n)-2)\log (-x)+2(\log{4\pi})m(i/n)-2(m(i/n))^2}{2\sqrt{m(i/n)}}I_0(x,y;m(i/n)) \nonumber
\\
&&
+\left(\frac{1}{n}\sum_{i=1}^{n}I_0(x,y;m(i/n)) +\frac{1}{n}\sum_{i=1}^{n} I_1(x,y;m(i/n))\right)o\left( \frac{1}{\log n} \right) \nonumber\\
&\sim&\frac{\log\log{n}}{2\sqrt{2\pi}\log{n}} \int_0^1 \sqrt{m(s)}\exp\left(-\frac{m(s)\lx{-}\log(xy)+\frac{(\log (x/y))^2}{4m(s)}}{2} \right)ds
\end{eqnarray}
as $n \rightarrow \infty$, where $I_0(x,y;m(i/n))$, $I_1(x,y;m(i/n))$ and $I_2(x,y;m(i/n))$ are given by \eqref{eq3.12}, \eqref{eq3.13} and \eqref{eq3.14}, respectively.

By Taylor expansion with Lagrange reminder term, we have
\begin{eqnarray}
\label{eq3.21}
&&
\Phi\left(\frac{\Phi^{-}(1+\frac{x}{n})-\rho_{ni}\Phi^{-}(1+\frac{t}{n})}{\sqrt{1-\rho_{ni}^2}}\right)\nonumber
\\
&=&\Phi(\sqrt{m(i/n)}+\frac{\log(t/x)}{2\sqrt{m(i/n)}})\nonumber
\\
&&
+\phi\left(\sqrt{m(i/n)}+\frac{\log(t/x)}{2\sqrt{m(i/n)}}\right) \left(\frac{\Phi^{-}(1+\frac{x}{n})-\rho_{ni}\Phi^{-}(1+\frac{t}{n})}{\sqrt{1-\rho_{ni}^2}}-\sqrt{m(i/n)}-\frac{\log(t/x)}{2\sqrt{m(i/n)}}\right)\nonumber
\\
&&
+\frac{1}{2}\theta_i \phi(\theta_i)
\left(\frac{\Phi^{-}(1+\frac{x}{n})-\rho_{ni}\Phi^{-}(1+\frac{t}{n})}{\sqrt{1-\rho_{ni}^2}}
-\sqrt{m(i/n)}-\frac{\log(t/x)}{2\sqrt{m(i/n)}}\right)^2,
\end{eqnarray}
where
\begin{eqnarray}
&&
\min\left(\frac{\Phi^{-}(1+\frac{x}{n})-\rho_{ni}\Phi^{-}(1+\frac{t}{n})}{\sqrt{1-\rho_{ni}^2}},
\sqrt{m(i/n)}+\frac{\log(t/x)}{2\sqrt{m(i/n)}}\right) \nonumber
\\
&<& \theta_i <
\max\left(\frac{\Phi^{-}(1+\frac{x}{n})-\rho_{ni}\Phi^{-}(1+\frac{t}{n})}{\sqrt{1-\rho_{ni}^2}},
\sqrt{m(i/n)}+\frac{\log(t/x)}{2\sqrt{m(i/n)}}\right). \nonumber
\end{eqnarray}

Note that
\begin{eqnarray}\label{addeq1}
&&
\frac{1}{n} \sum_{i=1}^{n} \int_y^{-\frac{1}{\log n}}
|\theta_i| \cdot \phi(\theta_i)
\left(\frac{\Phi^{-}(1+\frac{x}{n})-\rho_{ni}\Phi^{-}(1+\frac{t}{n})}{\sqrt{1-\rho_{ni}^2}}
-\sqrt{m(i/n)}-\frac{\log(t/x)}{2\sqrt{m(i/n)}}\right)^2dt\nonumber
\\
&\leq&
C \frac{1}{n} \sum_{i=1}^{n} \int_y^{-\frac{1}{\log n}}
\left( \frac{\log \log n}{4\log n} \left( \sqrt{m(i/n)}-\frac{\log(t/x)}{2\sqrt{m(i/n)}}\right) -\frac{(\log (-t))^2 +\left( \log 4\pi -3m(i/n) -2 \right)\log (-t)}{8\sqrt{m(i/n)}\log n}\right.\nonumber\\
&&\left.
+\frac{(\log (-x))^2 +\left( \log 4\pi +m(i/n) -2 \right)\log (-x)}{8\sqrt{m(i/n)}\log n} -\frac{(m(i/n))^{\frac{3}{2}}-(\log 4\pi)\sqrt{m(i/n)}}{4\log n}
-\frac{1+\log (-t)}{\log n} o(1)\right)^2dt
\nonumber
\\
&=&O\left(\left(\frac{\log \log n}{\log n}\right)^2\right).
\end{eqnarray}
Combining \eqref{eq3.20}, \eqref{eq3.21} and \eqref{addeq1}, we have
\begin{eqnarray}
\label{eq3.22}
&&
-\frac{1}{n} \sum_{i=1}^{n} \int_y^{-\frac{1}{\log n}}
\Big[\Phi\Big(\frac{\Phi^{-}(1+\frac{x}{n})-\rho_{ni}\Phi^{-}(1+\frac{t}{n})}{\sqrt{1-\rho_{ni}^2}}\Big)\nonumber
-\Phi\Big(\sqrt{m(i/n)}+\frac{\log(t/x)}{2\sqrt{m(i/n)}}\Big)\Big]dt\nonumber
\\
&=&-\frac{1}{n} \sum_{i=1}^{n} \int_y^{-\frac{1}{\log n}}
\Big[\phi\Big(\sqrt{m(i/n)}+\frac{\log(t/x)}{2\sqrt{m(i/n)}}\Big) \Big(\frac{\Phi^{-}(1+\frac{x}{n})-\rho_{ni}\Phi^{-}(1+\frac{t}{n})}{\sqrt{1-\rho_{ni}^2}}-\sqrt{m(i/n)}\nonumber
\\
&&
-\frac{\log(t/x)}{2\sqrt{m(i/n)}}\Big)
+\frac{1}{2}\theta_i \phi(\theta_i)
\Big(\frac{\Phi^{-}(1+\frac{x}{n})-\rho_{ni}\Phi^{-}(1+\frac{t}{n})}{\sqrt{1-\rho_{ni}^2}}
-\sqrt{m(i/n)}
-\frac{\log(t/x)}{2\sqrt{m(i/n)}}\Big)^2\Big]dt\nonumber
\\
&\sim&
\frac{\log\log{n}}{2\sqrt{2\pi}\log{n}} \int_0^1 \sqrt{m(s)}\exp\left(-\frac{m(s) \lx{-} \log(xy)+\frac{(\log (x/y))^2}{4m(s)}}{2} \right)ds
\end{eqnarray}
as $n\rightarrow \infty$.

Note that
\begin{equation}
\label{eq3.23}
\frac{1}{n} \sum_{i=1}^{n} \int_{-\frac{1}{\log n}}^0
\max\left(\Phi\Big(\frac{\Phi^{-}(1+\frac{x}{n})
-\rho_{ni}\Phi^{-}(1+\frac{t}{n})}{\sqrt{1-\rho_{ni}^2}}\Big), \Phi\Big(\sqrt{m(i/n)}+\frac{\log(t/x)}{2\sqrt{m(i/n)}}\Big)\right) dt
=O\Big(\frac{1}{\log n}\Big).
\end{equation}
Combining \eqref{eq3.22}, \eqref{eq3.23} with Lemma \ref{lem3.1}, we can get
\begin{eqnarray}
\label{eq3.24}
&&
\frac{1}{n} \sum_{i=1}^{n} \int_0^y
\Phi\Big(\frac{\Phi^{-}(1+\frac{x}{n})-\rho_{ni}\Phi^{-}(1+\frac{t}{n})}{\sqrt{1-\rho_{ni}^2}}\Big)dt
-\int_0^1 \int_0^y
\Phi\Big(\sqrt{m(s)}+\frac{\log(t/x)}{2\sqrt{m(s)}}\Big)dtds\nonumber
\\
&\sim&\frac{\log\log{n}}{2\sqrt{2\pi}\log{n}} \int_0^1 \sqrt{m(s)}\exp\left(-\frac{m(s)\lx{-} \log(xy)+\frac{(\log (x/y))^2}{4m(s)}}{2} \right)ds
\end{eqnarray}
as $n\rightarrow \infty$.

Since
\begin{eqnarray*}
& & x\int_{0}^{1}\Phi\left( \sqrt{m(s)} +\frac{\log (x/y)}{2\sqrt{m(s)}} \right) ds
+ y\int_{0}^{1}\Phi\left( \sqrt{m(s)} +\frac{\log (y/x)}{2\sqrt{m(s)}} \right) ds \nonumber\\
&=&
x+\int_{0}^{1} \int_{0}^{y}\Phi\left(\sqrt{m(s)} +\frac{\log (t/x)}{2\sqrt{m(s)}}  \right) dtds
\end{eqnarray*}
and from Liao et al.(2016) we have
\begin{eqnarray}\label{eq3.25}
&&
\P\Big(n\big(\max_{1\leq i \leq n}F_{1}(X_i)-1\big)\leq x,
n\big(\max_{1\leq i \leq n}F_{2}(Y_i)-1\big)\leq y \Big)-G(x,y)
\nonumber
\\
&=&G(x,y)\left( \sum_{i=1}^n \log{\P\left( F_1(X_i)\leq 1+\frac{x}{n},
F_2(Y_i)\leq 1+\frac{y}{n} \right)} -\log{G(x,y)} \right)(1+o(1))\nonumber\\
&=&G(x,y)\left( -\sum_{i=1}^n \left(  1-\P\left( F_1(X_i)\leq 1+\frac{x}{n},
F_2(Y_i)\leq 1+\frac{y}{n} \right) \right) \right.\nonumber\\
& & \left. -x\int_0^1 \Phi\left( \sqrt{m(s)}+\frac{\log(x/y)}{2\sqrt{m(s)}}\right) ds
-y\int_0^1 \Phi\left( \sqrt{m(s)}+\frac{\log(y/x)}{2\sqrt{m(s)}}\right)  ds  \right)(1+o(1))\nonumber
\\
&=&G(x,y)\left( -\sum_{i=1}^n \left( \P\left( F_1(X_i)> 1+\frac{x}{n} \right)  +\P\left( F_2(Y_i)>1+\frac{y}{n} \right) -\P\left( F_1(X_i)> 1+\frac{x}{n}, F_2(Y_i)>1+\frac{y}{n} \right) \right)
\right.\nonumber\\
&& \left.
-x -\int_{0}^{1} \int_{0}^{y}\Phi\left(\sqrt{m(s)} +\frac{\log (t/x)}{2\sqrt{m(s)}}  \right) dtds \right)(1+o(1))\nonumber
\\
&=&G(x,y)\left( \frac{1}{n}\sum_{i=1}^n \int_{0}^{y}
\Phi\left( \frac{\Phi^{-}\left( 1+\frac{x}{n}\right) -\rho_{ni}\Phi^{-}\left( 1+\frac{t}{n} \right)}{\sqrt{1-\rho^2_{ni}}} \right)dt
-\int_{0}^{1}\int_{0}^{y}\Phi\left( \sqrt{m(s)}+\frac{\log (t/x)}{2\sqrt{m(s)}} \right) dtds \right) \nonumber\\
&&
\times(1+o(1))\nonumber
\\
&\sim&\frac{\log\log{n}}{2\sqrt{2\pi}\log{n}} G(x,y) \int_0^1 \sqrt{m(s)}\exp\left(-\frac{m(s)\lx{-}\log(xy)+\frac{(\log (x/y))^2}{4m(s)}}{2} \right)ds
\end{eqnarray}
as $n\to \infty$, which complete the proof.
\qed

Next, using \eqref{eq3.11}-\eqref{eq3.15}, the proofs of results for cases: $\lim_{n\to \infty}\min_{1\leq i \leq n}m(i/n)=\infty$ and $\lim_{n\to \infty}\max_{1\leq i \leq n}m(i/n)=0$ are given
in the following.

\noindent
\textbf{Proof of Theorem \ref{the2.2}.}~~
For $y<t<-\frac{1}{\log n}$ with $x\leq0, y\leq0$, by using \eqref{eq3.16} and $\lim_{n\to \infty} \min_{1\leq i \leq n} m(i/n)=\infty$, we can get
\begin{eqnarray}
\label{eq3.26}
& & \frac{\Phi^{-}\left( 1+\frac{x}{n} \right) -\rho_{ni}\Phi^{-}\left( 1+\frac{t}{n} \right)}{\sqrt{1-\rho_{ni}^2}} \nonumber\\
&=& \sqrt{m(i/n)} + \frac{\log \frac{t}{x}}{2\sqrt{m(i/n)}}
+ \frac{\log\log n}{4\log n}\left( \frac{\log \frac{t}{x}}{2\sqrt{m(i/n)}} -\sqrt{m(i/n)} \right) \nonumber\\
& &
+\frac{(\log 4\pi -2)\log \frac{t}{x} + (\log (-t))^2 -(\log (-x))^2 +m(i/n)\log \frac{t}{x}}{8\sqrt{m(i/n)}\log n} \nonumber\\
& &
+\frac{(m(i/n))^{\frac{3}{2}}-2\sqrt{m(i/n)}\log (-t)}{4\log n}
+ o\left( \frac{(m(i/n))^{\frac{3}{2}}}{\log n} \right).
\end{eqnarray}

Note that the tail of normal distribution has the following expansion
\begin{equation}\label{eq3.27}
1-\Phi(x)= \frac{\phi(x)}{x}(1-x^{-2}+3x^{-4}-15x^{-6}+o(x^{-6}))
\end{equation}
as $x\to \infty$, c.f., Castro (1987). By using \eqref{eq3.27} and
$\lim_{n\to \infty}\frac{\max_{1\leq i \leq n}m(i/n)}{\log \log n}=0$, we have
\lx{\begin{eqnarray}\label{eq3.28}
&&\Phi\left( \frac{\Phi^{-}\left(1+\frac{x}{n}\right)-\rho_{ni}\Phi^{-}\left( 1-\frac{1}{n\log n} \right)}{\sqrt{1-\rho_{ni}^{2}}} \right)\nonumber\\
&=&1-\Phi\left(- \frac{\Phi^{-}\left(1+\frac{x}{n}\right)-\rho_{ni}\Phi^{-}\left( 1-\frac{1}{n\log n} \right)}{\sqrt{1-\rho_{ni}^{2}}} \right)\nonumber\\
&=&
-\frac{2\sqrt{m(i/n)}(-x\log n)^{\frac{1}{2}}}{\sqrt{2\pi}{\log \log n}}
\exp\left( -\frac{1}{2}m(i/n)-\frac{(\log(-x))\log \log n}{4m(i/n)}-\frac{(\log \log n)^2}{8m(i/n)} \right)
(1+o(1))\nonumber\\
\end{eqnarray}}
and
\begin{eqnarray}\label{eq3.29}
&&1-\Phi\left( \frac{\Phi^{-}\left(1+\frac{x}{n}\right)-\rho_{ni}\Phi^{-}\left( 1+\frac{y}{n} \right)}{\sqrt{1-\rho_{ni}^{2}}} \right)\nonumber\\
&=&
\frac{(x/y)^{\frac{1}{2}}}{\sqrt{2\pi}\sqrt{m(i/n)}}
\exp\left( -\frac{m(i/n)}{2}  \right) \left(1 -\frac{(\log{\frac{x}{y}})^2-4\log \frac{x}{y}+8}{8m(i/n)} +o\left( \frac{1}{m(i/n)} \right) \right)
\end{eqnarray}
for large $n$.

From \eqref{eq3.12}-\eqref{eq3.15}, $\lim_{n\to \infty} \min_{1\leq i \leq n} m(i/n)= \infty$
and $\lim_{n\to \infty} \frac{\max_{1\leq i\leq n}m(i/n)}{\log \log n}=0$,
it follows that
\begin{eqnarray*}
I_0(x,y;m(i/n))
&=& -2x\phi\left(  \sqrt{m(i/n)}+\frac{\log \frac{x}{y} }{2\sqrt{m(i/n)}} \right) \left( 1-\frac{\log{\frac{x}{y}}+2}{2m(i/n)}+o(\frac{1}{m(i/n)})  \right) ,
\end{eqnarray*}
\begin{eqnarray*}
I_1(x,y;m(i/n))
&=& 2x\phi\left(  \sqrt{m(i/n)}+\frac{\log \frac{x}{y} }{2\sqrt{m(i/n)}} \right) \left( \log(-y)-2+o(1) \right) ,
\end{eqnarray*}
\begin{eqnarray*}
I_2(x,y;m(i/n))
&=& -2x\phi\left(  \sqrt{m(i/n)}+\frac{\log \frac{x}{y} }{2\sqrt{m(i/n)}} \right) \left( (\log(-y))^2-4\log(-y)+8+o(1) \right) ,
\end{eqnarray*}
\begin{eqnarray*}
I_3(x,y;m(i/n))
&=& 2x\phi\left(  \sqrt{m(i/n)}+\frac{\log \frac{x}{y} }{2\sqrt{m(i/n)}} \right) \left( (\log(-y))^3 -6(\log(-y))^2 + 24 \log(-y) -48+o(1) \right),
\end{eqnarray*}
which implies that
\begin{eqnarray}
\label{eq3.30}
&&\int_{y}^{-\frac{1}{\log n}} t d \Phi\left( \frac{\Phi^{-}\left(1+\frac{x}{n}\right)-\rho_{ni}\Phi^{-}\left( 1+\frac{t}{n} \right)}{\sqrt{1-\rho_{ni}^2}} \right)\nonumber\\
&=&\int_{y}^{-\frac{1}{\log n}} \frac{(-t)\rho_{ni}}{n\sqrt{1-\rho_{ni}^2}}
\frac{\phi\left(  \frac{\Phi^{-}\left( 1+\frac{x}{n} \right)-\rho_{ni}\Phi^{-}\left( 1+\frac{t}{n} \right)}{\sqrt{1-\rho_{ni}^2}} \right)}{\phi\left(\Phi^{-}\left( 1+\frac{t}{n} \right)\right)}dt\nonumber\\
&=&\frac{1}{2\sqrt{m(i/n)}} \left(1-\frac{3m(i/n)}{4\log n} (1+o(1)) \right) \nonumber\\
& &\times \left[ \left(1+\frac{(1+m(i/n))\log \log n}{4\log n} -\frac{(m(i/n))^2}{4\log n}+o\left( \frac{(m(i/n))^2}{\log n} \right) \right)I_0(x,y;m(i/n)) \right.\nonumber\\
&&\left. +\left(\frac{1}{16\log n} -\frac{\log \log n+\log{4\pi}-2}{16m(i/n)\log n} \right) \left(I_2(x,y;m(i/n))+2\log (-x) I_1(x,y;m(i/n)) \right. \right.\nonumber\\
&&\left. \left.
+(\log(-x))^2 I_0(x,y;m(i/n))\right) +\frac{1}{16m(i/n)\log n} \left(I_3(x,y;m(i/n))+\log(-x)I_2(x,y;m(i/n))
\right. \right.\nonumber\\
&&\left. \left.
-(\log(-x))^2 I_1(x,y;m(i/n))-(\log(-x))^3 I_0(x,y;m(i/n))\right) \right.\nonumber\\
&&\left. -\left(\frac{m(i/n)}{4\log n}+o\left( \frac{m(i/n)}{4\log n} \right)  \right) I_1(x,y;m(i/n))
 \right]
\nonumber\\
&=&
\frac{(xy)^{\frac{1}{2}}}{\sqrt{2\pi}\sqrt{m(i/n)}}\exp\left( -\frac{m(i/n)}{2} \right)
\left( 1-\frac{\left( \log \frac{x}{y} \right)^2+ 4\log \frac{x}{y}+8}{8m(i/n)}
+o\left( \frac{1}{m(i/n)} \right)\right).
\end{eqnarray}
Hence, by using \eqref{eq3.28}-\eqref{eq3.30} we have
\begin{eqnarray}
\label{eq3.31}
&&\int_y^{-\frac{1}{\log n}} \left( 1-\Phi\left(  \frac{\Phi^{-}\left( 1+\frac{x}{n} \right)-\rho_{ni}\Phi^{-}\left( 1+\frac{t}{n} \right)}{\sqrt{1-\rho_{ni}^2}}  \right) \right)dt \nonumber\\
&=& -\frac{1}{\log n}\left( 1-\Phi\left(  \frac{\Phi^{-}\left(1+\frac{x}{n}\right)-\rho_{ni}\Phi^{-}\left( 1-\frac{1}{n\log n} \right)}{\sqrt{1-\rho_{ni}^2}} \right)  \right)\nonumber\\
&&-y\left(  1-\Phi\left(  \frac{\Phi^{-}\left(1+\frac{x}{n}\right)-\rho_{ni}\Phi^{-}\left( 1+\frac{y}{n} \right)}{\sqrt{1-\rho_{ni}^{2}}} \right) \right)\nonumber\\
&& + \int_{y}^{-\frac{1}{\log n}} t d \Phi\left( \frac{\Phi^{-}\left(1+\frac{x}{n}\right)-\rho_{ni}\Phi^{-}\left( 1+\frac{t}{n} \right)}{\sqrt{1-\rho_{ni}^2}}  \right)\nonumber\\
&=& \sqrt{\frac{2}{\pi}}\frac{(xy)^{\frac{1}{2}}}{\sqrt{m(i/n)}}\exp\left( -\frac{m(i/n)}{2} \right)(1+o(1))
\end{eqnarray}
as $n\to \infty$.

It follows from \eqref{eq3.23} and \eqref{eq3.31} that
\begin{eqnarray*}
&&\P(n(\max_{1\leq i \leq n}F_1(X_i)-1)\leq x,n(\max_{1\leq i \leq n}F_2(Y_i)-1)\leq y)
-e^{x+y}\nonumber\\
&=&e^{x+y}\left( \frac{1}{n} \sum_{i=1}^{n}\int_{y}^0 \left( 1-\Phi\left( \frac{\Phi^{-}\left( 1+\frac{x}{n} \right)-\rho_{ni}\Phi^{-}\left( 1+\frac{t}{n}
\right)}{\sqrt{1-\rho_{ni}^2}} \right) \right)dt \right)(1+o(1)) \nonumber\\
&=& e^{x+y}\left(  \frac{1}{n}\sum_{i=1}^{n}\frac{2(xy)^{\frac{1}{2}}}{\sqrt{2\pi}\sqrt{m(i/n)}}
\exp\left( -\frac{m(i/n)}{2} \right) \right)
(1+o(1)).
\end{eqnarray*}

The proof is complete.
\qed

\noindent
\textbf{Proof of Theorem \ref{the2.3}.}~~
Here we only prove the case of $x\neq y$ since the proof of case $x=y$ is similar. For $\max(x,y) \leq t \leq -\frac{1}{\log n}$, $x<0$, $y<0$, we have
\begin{eqnarray}\label{eq3.32}
& & \frac{\Phi^{-}\left( 1+\frac{\min(x,y)}{n} \right)-\rho_{ni}\Phi^{-}\left( 1+\frac{t}{n} \right)}{\sqrt{1-\rho_{ni}^2}} \nonumber\\
&=&\sqrt{m(i/n)}+\frac{\log \frac{t}{\min(x,y)}}{2\sqrt{m(i/n)}} -\frac{\log \frac{t}{\min(x,y)}}{4\sqrt{m(i/n)}\log n} +\frac{(\log 4\pi +\log \log n)\log \frac{t}{\min(x,y)}}{8\sqrt{m(i/n)}\log n} \nonumber\\
& & +\frac{(\log (-t))^2-(\log (-\min(x,y)))^2}{8\sqrt{m(i/n)}\log n} -\frac{\log \log n}{4\log n}\sqrt{m(i/n)} -\frac{\log (-t)}{2\log n}\sqrt{m(i/n)}\nonumber\\
& & +\frac{\sqrt{m(i/n)}\log \frac{t}{\min(x,y)}}{8\log n}+o\left( \frac{(\log \log n)\sqrt{m(i/n)}}{\log n} \right),
\end{eqnarray}
due to \eqref{eq3.16}, $\lim_{n\to \infty}\max_{1\leq i \leq n}m(i/n)=0$ and
$\lim_{n\to \infty} (\log \log n) \min_{1\leq i \leq n}m(i/n)=\infty$.
Noting that $\frac{\Phi^{-}\left( 1+\frac{\min(x,y)}{n}\right) -\rho_{ni}\Phi^{-}\left( 1+\frac{t}{n}  \right)}{\sqrt{1-\rho_{ni}^2}} \rightarrow -\infty$ for $t \in [\max(x,y),-\frac{1}{\log n}]$, we have
\begin{eqnarray}\label{eq3.33}
& & \Phi\left( \frac{\Phi^{-}\left( 1+\frac{\min(x,y)}{n}\right) -\rho_{ni}\Phi^{-}\left( 1-\frac{1}{n\log n}  \right)}{\sqrt{1-\rho_{ni}^2}} \right) \nonumber\\
&=& 1-\Phi\left( - \frac{\Phi^{-}\left( 1+\frac{\min(x,y)}{n}\right) -\rho_{ni}\Phi^{-}\left( 1-\frac{1}{n\log n}  \right)}{\sqrt{1-\rho_{ni}^2}} \right) \nonumber\\
&=&O\left( \frac{\sqrt{m(i/n)}(\log n)^{\frac{1}{2}}\exp\left( -\frac{\lx{(\log (-\min(x,y)\log n))^2}}{8m(i/n)} \right)}{\log \log n}  \right)
\end{eqnarray}
and
\begin{eqnarray}\label{eq3.34}
& & \Phi\left( \frac{\Phi^{-}\left( 1+\frac{\min(x,y)}{n}\right) -\rho_{ni}\Phi^{-}\left( 1+\frac{\max(x,y)}{n}  \right)}{\sqrt{1-\rho_{ni}^2}} \right) \nonumber\\
&=& 1-\Phi\left( - \frac{\Phi^{-}\left( 1+\frac{\min(x,y)}{n}\right) -\rho_{ni}\Phi^{-}\left( 1+\frac{\max(x,y)}{n}  \right)}{\sqrt{1-\rho_{ni}^2}} \right) \nonumber\\
&=& \frac{\sqrt{\frac{2}{\pi}}\sqrt{m(i/n)}\left( \frac{\min(x,y)}{\max(x,y)} \right)^{\frac{1}{2}}}{\left(\log \frac{\min(x,y)}{\max(x,y)}\right)} \exp\left(
- \frac{\left( \log \frac{\max(x,y)}{\min(x,y)}  \right)^2}{8m(i/n)}  \right)\nonumber\\
&& \times\left( 1-\frac{1}{2}m(i/n)+\frac{2m(i/n)}{\log \frac{\min(x,y)}{\max(x,y)}}-\frac{4m(i/n)}{\left( \log \frac{\max(x,y)}{\min(x,y)} \right)^2} +o(m(i/n)) \right).
\end{eqnarray}

From \eqref{eq3.12}-\eqref{eq3.15}, $\lim_{n\to \infty} \max_{1\leq i \leq n} m(i/n)=0$
and $\lim_{n\to \infty} (\log \log n) \min_{1\leq i\leq n}m(i/n)=\infty$,
it follows that
\begin{eqnarray*}
&&I_0(\min(x,y),\max(x,y);m(i/n))\nonumber\\
&=&-\frac{4\min(x,y)m(i/n)\phi\left(  \sqrt{m(i/n)}+\frac{\log \frac{\min(x,y)}{\max(x,y)} }{2\sqrt{m(i/n)}} \right)}{\left( \log \frac{\min(x,y)}{\max(x,y)} \right)\left(  1+\frac{2m(i/n)}{\log \frac{\min(x,y)}{\max(x,y)}} \right)}\left( 1-\frac{4m(i/n)}{\left( \log \frac{\min(x,y)}{\max(x,y)} \right)^2}(1+o(1))  \right)\nonumber\\
&&\times \left( 1-O\left( \frac{\exp\left( \lx{-\frac{(\log (-\min(x,y)\log n))^2}{8m(i/n)}+\frac{\left( \log \frac{\min(x,y)}{\max(x,y)}\right)^2}{8m(i/n)}} \right)}{(\log n)^{\frac{1}{2}}\log \log n} \right)  \right),
\end{eqnarray*}
\begin{eqnarray*}
&&I_1(\min(x,y),\max(x,y);m(i/n))\nonumber\\
&=&\frac{4\min(x,y)m(i/n)\phi\left(  \sqrt{m(i/n)}+\frac{\log \frac{\min(x,y)}{\max(x,y)} }{2\sqrt{m(i/n)}} \right)}{\left( \log \frac{\min(x,y)}{\max(x,y)} \right)\left(  1+\frac{2m(i/n)}{\log \frac{\min(x,y)}{\max(x,y)}} \right)}\left( 1-\frac{4m(i/n)}{\left( \log \frac{\min(x,y)}{\max(x,y)} \right)^2}(1+o(1))  \right)\nonumber\\
&&\times \left( \log(-\max(x,y))-\frac{4m(i/n)}{\log \frac{\min(x,y)}{\max(x,y)}} +o(m(i/n))  \right),
\end{eqnarray*}
\begin{eqnarray*}
&&I_2(\min(x,y),\max(x,y);m(i/n))\nonumber\\
&=&-\frac{4\min(x,y)m(i/n)\phi\left(  \sqrt{m(i/n)}+\frac{\log \frac{\min(x,y)}{\max(x,y)} }{2\sqrt{m(i/n)}} \right)}{\left( \log \frac{\min(x,y)}{\max(x,y)} \right)\left(  1+\frac{2m(i/n)}{\log \frac{\min(x,y)}{\max(x,y)}} \right)}\left( 1-\frac{4m(i/n)}{\left( \log \frac{\min(x,y)}{\max(x,y)} \right)^2}(1+o(1))  \right)\nonumber\\
&&\times \left((\log(-\max(x,y)))^2-\frac{8(\log (-\max(x,y)))m(i/n)}{\log \frac{\min(x,y)}{\max(x,y)}} +o(m(i/n))  \right)
\end{eqnarray*}
and
\begin{eqnarray*}
&&I_3(\min(x,y),\max(x,y);m(i/n))\nonumber\\
&=&\frac{4\min(x,y)m(i/n)\phi\left(  \sqrt{m(i/n)}+\frac{\log \frac{\min(x,y)}{\max(x,y)} }{2\sqrt{m(i/n)}} \right)}{\left( \log \frac{\min(x,y)}{\max(x,y)} \right)\left(  1+\frac{2m(i/n)}{\log \frac{\min(x,y)}{\max(x,y)}} \right)}\left( 1-\frac{4m(i/n)}{\left( \log \frac{\min(x,y)}{\max(x,y)} \right)^2}(1+o(1))  \right)\nonumber\\
&&\times \left((\log(-\max(x,y)))^3-\frac{12(\log (-\max(x,y)))^2m(i/n)}{\log \frac{\min(x,y)}{\max(x,y)}} +o(m(i/n))  \right),
\end{eqnarray*}
which implies that
\begin{eqnarray}\label{eq3.35}
&&\int_{\max(x,y)}^{-\frac{1}{\log n}} t d\Phi\left( \frac{\Phi^{-}\left( 1+\frac{\min(x,y)}{n} \right)-\rho_{ni}\Phi^{-}\left( 1+\frac{t}{n} \right)}{\sqrt{1-\rho_{ni}^2}}  \right)\nonumber\\
&=&\int_{\max(x,y)}^{-\frac{1}{\log n}} \frac{(-t)\rho_{ni}}{n\sqrt{1-\rho_{ni}^2}}
\frac{\phi\left(  \frac{\Phi^{-}\left( 1+\frac{x}{n} \right)-\rho_{ni}\Phi^{-}\left( 1+\frac{t}{n} \right)}{\sqrt{1-\rho_{ni}^2}} \right)}{\phi\left(\Phi^{-}\left( 1+\frac{t}{n} \right)\right)}dt\nonumber\\
&=&\frac{1-\frac{3m(i/n)}{4\log n}(1+o(1))}{2\sqrt{m(i/n)}} \nonumber\\
&& \times \left[
\left( 1+\frac{(1+m(i/n))\log \log n}{4\log n}+o\left( \frac{\log \log n}{\log n}m(i/n) \right)  \right)I_0\left( \min(x,y),\max(x,y);m(i/n)\right) \right. \nonumber\\
&&\left. + \left( \frac{2-\log 4\pi -\log \log n}{16m(i/n)\log n}+\frac{1}{16\log n} \right)\left( I_2(\min(x,y),\max(x,y);m(i/n)) \right.
\right. \nonumber\\
&&\left. \left. \quad +2(\log (-\min(x,y)))
I_1(\min(x,y),\max(x,y);m(i/n)) \right. \right. \nonumber\\
&&\left.\left. \quad + (\log (-\min(x,y)))^2I_0(\min(x,y),\max(x,y);m(i/n))\right) \right. \nonumber\\
&& \left. +\frac{1}{16m(i/n)\log n}\left( I_3(\min(x,y),\max(x,y);m(i/n))\right. \right. \nonumber\\
&&\left. \left. \quad
-(\log (-\min(x,y)))^2 I_1(\min(x,y),\max(x,y);m(i/n))
\right. \right. \nonumber\\
&&\left. \left. \quad + (\log (-\min(x,y))) I_2(\min(x,y),\max(x,y);m(i/n)) \right.\right. \nonumber\\
&&\left.\left. \quad- (\log (-\min(x,y)))^3 I_0(\min(x,y),\max(x,y);m(i/n)) \right)  \right. \nonumber\\
&&\left. +o\left(\frac{\log \log n}{\log n}\right) I_1(\min(x,y),\max(x,y);m(i/n)) \right] \nonumber\\
&=&\frac{\sqrt{\frac{2}{\pi}}\sqrt{m(i/n)}(xy)^{\frac{1}{2}}}{\left( \log \frac{\min(x,y)}{\max(x,y)} \right) }\exp\left( -\frac{\left( \log{\frac{\min(x,y)}{\max(x,y)}} \right)^2}{8m(i/n)} \right)\nonumber\\
&& \times
\left( 1-\frac{1}{2}m(i/n)-\frac{2m(i/n)}{\log{\frac{\min(x,y)}{\max(x,y)}}} -\frac{4m(i/n)}{\left(\log{\frac{\min(x,y)}{\max(x,y)}} \right)^2} +o(m(i/n)) \right).
\end{eqnarray}
Hence, it follows from \eqref{eq3.33}-\eqref{eq3.35} that
\begin{eqnarray}\label{case3}
&&\int_{\max(x,y)}^{-\frac{1}{\log n}}  \Phi\left( \frac{\Phi^{-}\left( 1+\frac{\min(x,y)}{n} \right)-\rho_{ni}\Phi^{-}\left(1+\frac{t}{n}  \right)}{\sqrt{1-\rho_{ni}^2}}  \right)dt \nonumber\\
&=&-\frac{1}{\log n} \Phi\left( \frac{\Phi^{-}\left( 1+\frac{\min(x,y)}{n} \right)-\rho_{ni}\Phi^{-}\left(1-\frac{1}{n\log n}  \right)}{\sqrt{1-\rho_{ni}^2}}  \right)  \nonumber\\
&&
- \max(x,y)\Phi\left( \frac{\Phi^{-}\left( 1+\frac{\min(x,y)}{n} \right)-\rho_{ni}\Phi^{-}\left(1+\frac{\max(x,y)}{n}  \right)}{\sqrt{1-\rho_{ni}^2}}  \right) \nonumber\\
&&
-\int_{\max(x,y)}^{-\frac{1}{\log n}} t d \Phi\left( \frac{\Phi^{-}\left( 1+\frac{\min(x,y)}{n} \right)-\rho_{ni}\Phi^{-}\left(1+\frac{t}{n}  \right)}{\sqrt{1-\rho_{ni}^2}}  \right) \nonumber\\
&=&
\frac{8{m(i/n)}^{\frac{3}{2}}(xy)^{\frac{1}{2}}}{\sqrt{2\pi} \log{\frac{\min(x,y)}{\max(x,y)}}} \exp\left( -\frac{\left( \log{\frac{\min(x,y)}{\max(x,y)}} \right)^2}{8m(i/n)} \right)(1+o(1)).
\end{eqnarray}
Hence, by using \eqref{case3} and \eqref{eq3.23}, we have
\begin{eqnarray*}
&&\P\Big(n(\max_{1\leq i \leq n}F_1(X_i)-1)\leq x, n(\min_{1\leq i\leq n}F_2(Y_i)-1)\leq y\Big)-e^{\min(x,y)} \nonumber\\
&=&e^{\min(x,y)}\left( -\frac{1}{n}\sum_{i=1}^n\int_{\max(x,y)}^{-\frac{1}{\log n}}
\Phi\left( \frac{\Phi^{-}\left( 1+\frac{\min(x,y)}{n} \right)-\rho_{ni}\Phi^{-}\left(1+\frac{t}{n}  \right)}{\sqrt{1-\rho_{ni}^2}} \right)dt +O\left( \frac{1}{\log n} \right)\right)(1+o(1))\nonumber\\
&=&
-e^{\min(x,y)} \left( \frac{1}{n} \sum_{i=1}^{n} \frac{8\left(m(i/n)\right)^{\frac{3}{2}} (xy)^{\frac{1}{2}}}{\sqrt{2\pi} \log{\frac{\min(x,y)}{\max(x,y)}}} \exp\left( -\frac{\left( \log{\frac{\min(x,y)}{\max(x,y)}} \right)^2}{8m(i/n)} \right)  \right)(1+o(1)),
\end{eqnarray*}
which is the desired result. The proof is complete.
\qed

\noindent
{\bf Acknowledgements}~~This work was supported by the National Natural Science Foundation of China (grant No. 11601330), and the Funding Program for Junior
Faculties of College and Universities of Shanghai Education Committee (grant No. ZZslg16020).

\end{document}